\documentclass[a4paper,10pt, oneside]{article}
\usepackage{amsmath, amssymb, amsthm, color}
\usepackage[font=small,labelfont=md,textfont=it]{caption}
\usepackage{booktabs,makecell,multirow}
\usepackage{floatrow}
\usepackage[titletoc, title]{appendix}
\usepackage[colorlinks,linkcolor=blue,citecolor=blue]{hyperref}
\usepackage[capitalise, nosort]{cleveref}

\crefname{equation}{}{}
\crefname{lem}{Lemma}{Lemmas}
\crefname{thm}{Theorem}{Theorems}
\numberwithin{equation}{section}

\newcommand{\dual}[1]{\left\langle #1 \right\rangle}
\newcommand{\nm}[1]{\left\Vert #1 \right\Vert}
\newcommand{\snm}[1]{\left\vert #1 \right\vert}
\newcommand{\ssnm}[1]
{
  \left\vert\kern-0.25ex
    \left\vert\kern-0.25ex
      \left\vert
        #1
      \right\vert\kern-0.25ex
    \right\vert\kern-0.25ex
  \right\vert
}
\newcommand{\jmp}[1]{[\kern-.15em[ #1 ]\kern-.15em]}

\newtheorem{corollary}{Corollary}[section]
\newtheorem{Def}{Definition}[section]
\newtheorem{lem}{Lemma}[section]
\newtheorem{rem}{Remark}[section]
\newtheorem{thm}{Theorem}[section]

\begin{document}

\title{ \Large\bf Analysis of a time-stepping discontinuous Galerkin method for modified
anomalous subdiffusion problems
\thanks
  {
    This work was supported in part by National Natural Science Foundation
    of China (11771312).
  }
}
\author{
	Binjie Li \thanks{Email: libinjie@scu.edu.cn},
	Hao Luo \thanks{Corresponding author. Email: galeolev@foxmail.com},
	Xiaoping Xie \thanks{Email: xpxie@scu.edu.cn} \\
	{School of Mathematics, Sichuan University, Chengdu 610064, China}
}

\date{}
\maketitle

\begin{abstract}
  This paper analyzes a time-stepping discontinuous Galerkin method for 
  modified anomalous subdiffusion problems with two time fractional derivatives
  of orders $ \alpha $ and $ \beta $ ($ 0 < \alpha < \beta < 1 $). The stability
  of this method is established, the temporal accuracy of $ O(\tau^{m+1-\beta/2}) $
  is derived, where $m$ denotes the degree of polynomials for the temporal discretization. It is shown that, even the solution has singularity near $ t
  = {0+} $, this temporal accuracy can still be achieved by using the graded
  temporal grids. Numerical experiments are performed to verify
  the theoretical results.
\end{abstract}


\section{Introduction}
This paper considers the following modified anomalous fractional subdiffusion
problem:
\begin{equation}
  \label{eq:model}
  \left\{
    \begin{aligned}
      \partial_t u -
      \left( \kappa_1 D_{0+}^\alpha + \kappa_2 D_{0+}^\beta \right)
      \partial_x^2 u & = f   &  & \text{in $ \Omega \times (0,T) $,} \\
      u              & = 0   &  & \text{on $ \partial\Omega $,} \\
      u(\cdot,0)     & = u_0 &  & \text{in $\Omega$,}
    \end{aligned}
  \right.
\end{equation}
where $ \Omega \subset \mathbb R $ is an open
interval, $ \kappa_1 $ and $ \kappa_2 $ are two positive constants, $ 0 < T < \infty $, $ 0 < \alpha < \beta < 1 $, $ u_0 \in H_0^1(\Omega) $,
and $ f \in L^1(0,T;H^{-1}(\Omega)) $.

A considerable amount of research has been devoted to the numerical treatment
of time fractional diffusion problems, especially in the past decade. So far,
most of the existing algorithms are classified as fractional difference methods,
since they employ the $ L^1 $ formula, Gr\"unwald-Letnikov discretization or
fractional linear multi-step method to discretize the fractional derivatives.
Despite their ease of implementation, the fractional difference methods are
generally of temporal accuracy orders not greater than two; see \cite{Yuste2005,
Langlands2005, Yuste2006, Chen2007, Lin2007, Zhuang2008, Deng2009, Zhang2009524,
Chen2010,Cui2009, Li2011Numerical, Liu2011,Gao2011, LiXu2012,Zeng2013, Jin2014,
Wang2014, gao2014new, Li2016} and the references therein. We also note that Gao et.~al \cite{gao2014new} designed a formula to approximate the Caputo fractional
derivative of order $ \alpha $ ($ 0<\alpha<1 $) and applied this formula to
numerically solve time fractional diffusion problems; however, the theory of
stability and convergence was not established there. To improve the temporal accuracy,
fractional spectral methods, namely those algorithms using spectral methods to
discretize the fractional derivatives, were proposed; see \cite{Li2009,
Zheng2015,Li2017A}. Recently, Mustapha and Mclean (\cite{Mclean2009Convergence,
Mustapha2011Piecewise, Mustapha2015Time}) used the discontinuous Galerkin method
to approximate the time fractional derivatives, and they proposed a class of
methods that possess high-order temporal accuracy. Moreover, as the fractional
difference methods, the numerical solutions of their methods are computed in a
step by step fashion.

Due to the nonlocal property of the fractional derivatives, the computation and
storage cost of an accurate numerical solution to a time fractional diffusion
problem significantly exceeds that to a corresponding normal diffusion problem.
Naturally, developing high-order accuracy algorithms, especially those with
high-order temporal accuracy, is an efficient way to reduce the cost. However, as
aforementioned, generally the best temporal accuracy order of the fractional
difference methods is merely two. This motivates us to develop algorithms that
possess high-order accuracy in both space and time while retaining the advantage
of the fractional finite difference methods.

Following the work of \cite{Mustapha2015Time} for fractional diffusion equations, we analyze a time-stepping
discontinuous Galerkin method for problem \cref{eq:model}. Firstly, we establish
a new stability estimate. Secondly, we prove that the temporal accuracy is $
\mathcal O(\tau^{m+1-\beta/2}) $, and that if $ u $ has singularity near $ t =
{0+} $, then this temporal accuracy can still be achieved by using  graded
temporal grids.  We note that on  appropriate graded temporal grids, 
 \cite{Mustapha2015Time}
obtained the temporal accuracy $ O(\tau^{2-\beta/2}) $ in the case of $ m = 1 $
and  the temporal accuracy $ O(\tau^{m+(1-\beta)/2}) $ in the case of $ m
\geqslant 2 $.

The rest of this paper is organized as follows. \cref{sec:preli} introduces 
some  notations. 
\cref{sec:main,sec:proof}
establish the stability and convergence of the time-stepping discontinuous
Galerkin method. \cref{sec:numer} performs several numerical experiments to verify
the theoretical results. Finally, \cref{sec:concluding} gives 
concluding remarks.


\section{ Notation}
\label{sec:preli}
Let us first introduce the Riemann-Liouville fractional calculus operators.
\begin{Def}
  For $ 0 < \gamma < \infty $ and  any $ v \in L^1(0,T;X) $,  define
  \begin{align*}
    \left(I_{0+}^{\gamma,X} v\right)(t) &:=
    \frac1{ \Gamma(\gamma) }
    \int_0^t (t-s)^{\gamma-1} v(s) \, \mathrm{d}s, \quad 0 < t < T, \\
    \left(I_{T-}^{\gamma,X} v\right)(t) &:=
    \frac1{ \Gamma(\gamma) }
    \int_t^T (s-t)^{\gamma-1} v(s) \, \mathrm{d}s, \quad 0 < t < T,
  \end{align*}
where $ \Gamma(\cdot) $ is the gamma function.
\end{Def}
\begin{Def}
  For $ 0 < \gamma < 1 $, define
  \begin{align*}
    D_{0+}^{\gamma,X} & := D I_{0+}^{1-\gamma,X}, \\
    D_{T-}^{\gamma,X} & := -D I_{T-}^{1-\gamma,X},
  \end{align*}
  where $ D $ is the first-order differential operator in the distribution
  sense.
\end{Def}
\noindent Above $ X $ is a Banach space and $ L^1(0,T;X) $ is a standard $ X
$-valued Bochner $ L^1 $ space. For convenience, we shall simply use $
I_{0+}^\gamma $, $ I_{T-}^\gamma $, $ D_{0+}^\gamma $ and $ D_{T-}^\gamma $,
without indicating the underlying Banach space $ X $.

Next we introduce some vector valued spaces. Let $ X $ be a separable Hilbert
space with an inner product $ (\cdot,\cdot)_X $ and an orthonormal basis $ \{
e_j:\ j \in \mathbb N \} $, and let $ \mathcal O $ be an interval. For $ 0 <
\gamma < \infty $, define
\[
  H^\gamma(\mathcal O;X) := \left\{
    v \in L^2(\mathcal O;X):\
    \sum_{j=0}^\infty \nm{(v,e_j)_X}_{H^\gamma(\mathcal O)}^2 < \infty
  \right\}
\]
and equip this space with the norm
\[
  \nm{\cdot}_{H^\gamma(\mathcal O; X)} := \left(
    \sum_{j=0}^\infty \nm{(\cdot,e_j)_X}_{H^\gamma(\mathcal O)}^2
  \right)^\frac12,
\]
where $ L^2(\mathcal O;X) $ is an $ X $-valued Bochner $L^2$ space. If $ 0 <
\gamma < 1/2 $, we also introduce the   seminorm
\[
  \snm{\cdot}_{H^\gamma(\mathcal O;X)} := \left(
    \sum_{j=0}^\infty \snm{(\cdot,e_j)_X}_{H^\gamma(\mathcal O)}^2
  \right)^\frac12.
\]
Here, $ H^\gamma(\mathcal O) $ is a standard Sobolev space (see
\cite{Tartar2007}), and
\[
  \snm{v}_{H^\gamma(\mathcal O)} := \left(
    \int_\mathbb R \snm{\xi}^{2\gamma}
    \snm{\mathcal F(v\chi_{\mathcal O})(\xi)}^2 \,\mathrm{d}\xi
  \right)^\frac12
\]
for each $ v \in H^\gamma(\mathcal O) $ with $ 0 < \gamma < 1/2 $, where $
\mathcal F: L^2(\mathbb R) \to L^2(\mathbb R) $ is the Fourier transform
operator and $ \chi_{\mathcal O} $ is the indicator function of the interval $
\mathcal O $. For $ v \in H^i(\mathcal O;X) $ with $ i \in \mathbb N_{>0} $,
define its $ i $th weak derivative $ v^{(i)} $ by
\[
  v^{(i)} := \sum_{j=0}^\infty c_j^{(i)}(t) e_j, \quad t \in \mathcal O,
\]
where $ c_j := (v, e_j)_X $ and $ c_j^{(i)} $ is its $ i $th weak derivative. In
particular, $ v^{(1)} $ is abbreviated to $ v' $.

Additionally, for $ 0 \leqslant \delta < 1 $, define
\[
  L_\delta^2(\mathcal O;X) := \left\{
    v \in L^1(\mathcal O;X) :\
    \nm{v}_{L_\delta^2(\mathcal O;X)} < \infty
  \right\},
\]
where
\[
  \nm{v}_{L_\delta^2(\mathcal O;X)} := \left(
    \int_\mathcal O \snm{t}^\delta \nm{v(t)}_X^2 \, \mathrm{d}t
  \right)^\frac12.
\]
Conventionally, $ C(\mathcal O;X) $ is the set of all $ X $-valued continuous
functions defined on $ \mathcal O $, and $ P_j(\mathcal O;X) $ is the set of all
$ X $-valued polynomials of degree $ \leqslant j $  defined on $ \mathcal O $.
For convenience, $ \nm{\cdot}_{L_\delta^2(\mathcal O;\mathbb R)} $ and $
P_j(\mathcal O;\mathbb R) $ are abbreviated to $ \nm{\cdot}_{L_\delta^2(\mathcal
O)} $ and $ P_j(\mathcal O) $, respectively.

Now we introduce a temporal and spatial discretization space. Let $ \mathcal
K_h $ be a partition of $ \Omega $ consisting of open intervals, and let $ h
$  denote the maximum length of the elements in $ \mathcal K_h $. We introduce a graded mesh subdivision of  the temporal interval $(0,T)$. For $ J \in
\mathbb N_{>0} $ and $ \sigma \geqslant 1 $, we set 
\begin{equation}\label{graded}
\left\{
\begin{array}{ll}
 t_j :=
(j/J)^\sigma T \quad & \text{for   } 0 \leqslant j \leqslant J ,\\
 I_j :=
(t_{j-1},t_j) ,\ \tau_j := t_j - t_{j-1} \quad & \text{for   }   1 \leqslant j
\leqslant J ,
\end{array}
\right.
\end{equation}
 and use $ \tau $ to abbreviate $ \tau_J $. 
 Define
\begin{align*}
  S_h &:= \left\{
    v_h \in H_0^1(\Omega):\
    v_h|_K \in P_n(K),\ \forall K \in \mathcal K_h
  \right\}, \\
  W_{h,\tau} &:= \left\{
    V \in L^2(0,T;S_h):\
    V|_{I_j} \in P_m(I_j;S_h),\ \forall 1 \leqslant j \leqslant J
  \right\},
\end{align*}
where $ m \in \mathbb N $ and $ n \in \mathbb N_{>0} $. Moreover, for $ V \in
W_{h,\tau} $ we introduce the following notation:
\begin{align*}
  & V_j^{+}  := \lim_{t \to {t_j}+} V(t)
  \quad\text{for } 0 \leqslant j < J; \\
  & V_j  := \lim_{t \to {t_j}-} V({t}) =: V(t_j)
  \quad\text{for } 0 < j \leqslant J \text{, and }
  V_0 := 0; \\
  & \jmp{V_j}  := V_j^{+} - V_j
  \quad\text{for } 0 \leqslant j < J.
\end{align*}
We note that in the context, $ H^s(\Omega) $ ($ s \in \mathbb Z $) and $ H_0^s(\Omega)
$ ($ s \in \mathbb N_{>0} $) denote two standard Sobolev spaces (see
\cite{Tartar2007}).

Throughout this paper, we make the following conventions: each $ v \in L^1(\Omega\times(0,T))
$ is regarded as an element of $ L^1(0,T;L^1(\Omega)) $, also denoted by $ v $;
the notation $ a \lesssim b $ means that there exists a positive constant $ C $
depending only on $ \alpha $, $ \beta $, $ m $ or $ n $ such that $ a \leqslant
C b $, and $ a \sim b $ means $ a \lesssim b \lesssim a $; by $ C_{\times} $ we
denote a positive constant that only depends on $ {\times} $ and its value may
differ at each of its occurrences; if $ \mathcal O $ is a Lebesgue measurable
set of $ \mathbb R $ or $ \mathbb R^2 $, then $ \dual{v,w}_\mathcal O $ means $
\int_\mathcal O vw $; if $ X $ is a Banach space, then $ \dual{\cdot,\cdot}_X $
denotes the duality pairing between $ X' $ and $ X $.


\section{Main Results}
\label{sec:main}
Let us first describe the time-stepping discontinuous Galerkin method to be
analyzed as follows: seek $ U \in W_{h,\tau} $ such that
\begin{equation}
\label{eq:U}
\begin{aligned}
&
\dual{U', V}_{\Omega_T} +
\sum_{j=0}^{J-1} \dual{\jmp{U_j}, V_j^{+}}_\Omega +
\dual{
    \left(
    \kappa_1 D_{0+}^\alpha +
    \kappa_2 D_{0+}^\beta
    \right) \partial_x U, \partial_x V
}_{\Omega_T} \\
=&
\dual{R_hu_0, V_0^{+}}_\Omega +
\dual{f,V}_{L^\infty(0,T;H_0^1(\Omega))}
\end{aligned}
\end{equation}
for all $ V \in W_{h,\tau} $, where $ \Omega_T := \Omega \times (0,T) $ and the
projection operator $ R_h $ is defined by
\[
\dual{\partial_x(v-R_hv), \partial_x v_h}_{\Omega} = 0,
\quad \forall v \in H_0^1(\Omega),\ \forall v_h \in V_h.
\]
Above $ U' $ is understood by
\[
U'|_{I_j} := \left( U|_{I_j} \right)',
\quad 1 \leqslant j \leqslant J.
\]

Then, assuming $ X $ to be a Banach space, we define an interpolation operator $
Q_\tau^X $ as follows \cite[Chapter~12]{Thomee2006}: given $ v \in C((0,T];X)
\cap L^1(0,T;X) $, the interpolant $ Q_\tau^X v $ fulfills, for each $ 1 \leqslant j \leqslant J $,
\[
  \left\{
    \begin{aligned}
      & \left( Q_\tau^X v \right)|_{I_j} \in P_m(I_j;X),
      \quad \lim_{t\to t_{j-}}\left( Q_\tau^Xv \right)(t) = v(t_j), \\
      & \int_{t_{j-1}}^{t_j} \left( v-Q_\tau^Xv \right) q \, \mathrm{d}t = 0
      \quad \text{for all } q \in P_{m-1}(I_j).
    \end{aligned}
  \right.
\]
 Below we will use $ Q_\tau $ instead of
$ Q_\tau^X $ when no confusion will arise.

Now we are ready to state the main results of this paper, and, for convenience,
we assume that $ u $ is the solution to problem \cref{eq:model}. 

\begin{thm}
  \label{thm:stability} The scheme \eqref{eq:U} admits a unique solution $U$. In addition,
 if $ f \in L_\beta^2(0,T;H^{-1}(\Omega)) $, then
  \begin{equation}
    \label{eq:stab}
    \begin{aligned}
      & \nm{U(t_j)}_{L^2(\Omega)} +
      \sqrt{\kappa_1} \snm{U}_{H^{\alpha/2}(0,t_j;H_0^1(\Omega))} +
      \sqrt{\kappa_2} \snm{U}_{H^{\beta/2}(0,t_j;H_0^1(\Omega))} \\
      \lesssim{} &
      \nm{u_0}_{H_0^1(\Omega)} + 1/\sqrt{\kappa_2}
      \nm{f}_{L_\beta^2(0,t_j;H^{-1}(\Omega))}
    \end{aligned}
  \end{equation}
  for each $ 1 \leqslant j \leqslant J $.
\end{thm}

\begin{thm}
  \label{thm:conv}
  If $ u' \in L_\beta^2(0,T;H_0^1(\Omega)) $, then
  \begin{equation}
    \label{eq:conv}
    \begin{aligned}
      {} &
      \nm{\theta(t_j)}_{L^2(\Omega)} + \sqrt{\kappa_1}
      \snm{\theta}_{H^{\alpha/2}(0,t_j;H_0^1(\Omega))} +
      \sqrt{\kappa_2}
      \snm{\theta}_{H^{\beta/2}(0,t_j;H_0^1(\Omega))} \\
      \lesssim{} &
      \eta_{j,1} + \eta_{j,2} + \eta_{j,3}
    \end{aligned}
  \end{equation}
  for each $ 1 \leqslant j \leqslant J $, where $ \theta := U - Q_\tau R_h u $
  and
  \begin{align*}
    \eta_{j,1} &:= h^{\min\{2,n\}}/\sqrt{\kappa_2}
    \nm{(I-R_h)u'}_{L_\beta^2(0,t_j;H_0^1(\Omega))}, \\
    \eta_{j,2} &:= \sqrt{\kappa_1} \left(
      \sum_{i=1}^j \tau_i^{2-\alpha}
      \inf_{0\leqslant\delta<1} \frac{t_i^{-\delta}}{1-\delta}
      \nm{(u-Q_\tau u)'}^2_{L^2_\delta(I_i;H_0^1(\Omega))}
    \right)^\frac12, \\
    \eta_{j,3} &:= \sqrt{\kappa_2} \left(
      \sum_{i=1}^j \tau_i^{2-\beta}\inf_{0\leqslant\delta<1}
      \frac{t_i^{-\delta}}{1-\delta}
      \nm{(u-Q_\tau u)'}^2_{L^2_\delta(I_i;H_0^1(\Omega))}
    \right)^\frac12.
  \end{align*}
\end{thm}

\begin{corollary}
  \label{coro:conv}
  If $ u \in H^{m+1}(0,T;H^1(\Omega)) $ and
  \[
    u' \in L_\beta^2(0,T;H_0^1(\Omega) \cap H^{n+1}(\Omega)),
  \]
  then
  \[
    \nm{(u-U)(t_j)}_{L^2(\Omega)} \lesssim \nu_{j,1} + \nu_{j,2}
  \]
  for each $ 1 \leqslant j \leqslant J $, where
  \begin{align*}
    \nu_{j,1} &:= h^{\min\{2,n\}+n} / \sqrt{\kappa_2}
    \nm{u'}_{L_\beta^2(0,t_j;H^{n+1}(\Omega))} +
    h^{n+1} \nm{u(t_j)}_{H^{n+1}(\Omega)}, \\
    \nu_{j,2} &:=
    \left(
      \sqrt{\kappa_1} \tau_j^{1+m-\alpha/2} +
      \sqrt{\kappa_2} \tau_j^{1+m-\beta/2}
    \right)
    \nm{u}_{H^{m+1}(0,t_j;H^1(\Omega))}.
  \end{align*}
\end{corollary}

\begin{corollary}
  \label{coro:graded}
  If $ u(x, t) = t^r \phi(x) $ for  $ (x,t)\in \Omega_T$, with $ (\beta-1)/2 < r \leqslant m+1/2 $
  and $ \phi \in H_0^1(\Omega) \cap H^{n+1}(\Omega) $, then
  \[
    \begin{aligned}
      \nm{(u - U)(t_j)}_{L^2(\Omega)}
      \lesssim {}& C_{\sigma,r}\left(
        \sqrt{\kappa_1}  +
        \sqrt{\kappa_2}
      \right)\epsilon_j \nm{\phi}_{H_0^1(\Omega)}+h^{n+1} \,t_j^r \nm{\phi}_{H^{n+1}(\Omega)} \\
      {}&\quad +
      h^{\min\{2,n\}+n} \,t_j^{r+(\beta-1)/2} \nm{\phi}_{H^{n+1}(\Omega)} 
    \end{aligned}
  \]
  for all $ 1 \leqslant j \leqslant J $, where
  \[
    \epsilon_j:=
    \begin{cases}
      T^{(1-\sigma)(r+(1-\beta)/2)}\tau^{\sigma(r+(1-\beta)/2)}&
      \text{ if } \sigma < \sigma^*, \\
      (1+\ln\left(j\right)) T^{r-m-1/2}\tau^{m+1-\beta/2}&
      \text{ if } \sigma = \sigma^*, \\
      T^{r-m-1/2}\tau^{m+1-\beta/2}
      &
      \text{ if } \sigma >\sigma^*,
    \end{cases}
  \]
  $\sigma$ is the graded parameter in \eqref{graded},
and 
  \begin{equation}
    \label{eq:s1}
    \sigma^* := \frac{2m+2-\beta}{2r+1-\beta}\geqslant 1 .
  \end{equation}
\end{corollary}

\begin{rem}
    Due to the fact that
    \begin{align*}
        \nm{U}_{H^{\beta/2}(0,T;H_0^1(\Omega))} \lesssim C_T
        \snm{U}_{H^{\beta/2}(0,T;H_0^1(\Omega))},
    \end{align*}
    by \cref{thm:stability,thm:conv} we can also derive  the stability and
    error estimates of $ U $ with respect to the norm on $
    H^{\beta/2}(0,T;H_0^1(\Omega)) $.
\end{rem}

\begin{rem}
    If $ n \geqslant 2 $ and the condition of \cref{coro:conv} is satisfied, then
    \cref{thm:conv} implies
    \[
    \nm{ (U-Q_\tau R_hu)(t_j) }_{L^2(\Omega)} =
    O( h^{n+2} + \tau^{m+1-\beta/2} ).
    \]
    Assume that $ \mathcal K_h $ is quasi-uniform and $ \{x_i: 1 \leqslant i
    \leqslant N \} $ is the set of all nodes of $ \mathcal K_h $. Using the
    standard result
    \[
    R_hu(x_i,t_j) = u(x_i,t_j), \quad 1 \leqslant i \leqslant N,
    \]
    we obtain
    \[
    \max_{1\leqslant i\leqslant N} \snm{U(x_i,t_j)-u(x_i,t_j)} =
    O( h^{n+1} + h^{-1}\tau^{m+1-\beta/2} ).
    \]
    Therefore, if $ \tau $ is sufficiently small, then
    \[
    \max_{1\leqslant i\leqslant N}  \snm{U(x_i,t_j)-u(x_i,t_j)} = O(h^{n+1}).
    \]
\end{rem}
\begin{rem}
Though the graded grids are assumed in  \eqref{graded}, from the proofs in \cref{sec:proof} it is easy to see  that  \cref{thm:stability}, \cref{thm:conv}, and \cref{coro:conv} still hold for  more general temporal grids, with $\tau_j$ in  \cref{coro:conv}  replaced by $\max\limits_{1\leqslant i \leqslant j}\tau_i$.
\end{rem}
\begin{rem}
   First, \cref{coro:graded} shows that if $ u $ has singularity near $ t = {0+} $, then
    the graded grids in the time direction can improve the temporal accuracy to $ O(\tau^{m+1-\beta/2}) $ up to an factor $\ln(j)$ provided that $ \sigma = \sigma^* $. Numerical results show that our estimates are sharp for $ \sigma \leqslant \sigma^* $. Second, theoretically we can not expect the optimal accuracy $ O(\tau^{m+1}) $ as $ \sigma $ increases. However, numerical tests indicate that the optimal convergence rate can also be obtained if
      \begin{equation}\label{eq:s2}
       \sigma = \sigma^{**} := \frac{2m+2}{2r+1-\beta}>\sigma^*.
       \end{equation}
\end{rem}

\begin{rem}
We note that for the time stepping discontinuous Galerkin discretization of fractional diffusion problems,   \cite{Mustapha2015Time} obtained the temporal accuracy $ O(\tau^{2-\beta/2}) $ for $ m = 1 $ and    $ O(\tau^{m+(1-\beta)/2}) $  for $ m
\geqslant 2 $ on appropriate graded temporal grids.
\end{rem}

The rest of this section will briefly discuss the singularity of the solution to
problem \cref{eq:model} near $ t = {0+} $. Let $ \{\phi_j\}^\infty_{j = 0} $ be
an orthonormal basis of $ L^2(\Omega) $ such that $ \phi_j \in H_0^1(\Omega) $ and
\[
- \partial_{x}^2\phi_j  = \lambda_j \phi_j \quad x\in \Omega,
\]
where $ \{ \lambda_j \}_{j=0}^\infty \subset \mathbb R_{>0} $ is an
non-decreasing sequence. For each $ j \in \mathbb N $, define
\begin{align*}
  u_j(t) &:= \dual{u(x,t),\phi_j}_\Omega, \quad 0 < t < T, \\
  f_j(t) &:= \dual{f(x,t),\phi_j}_\Omega, \quad 0 < t < T.
\end{align*}
Evidently, $ u_j $ satisfies the fractional ordinary equation
\begin{equation}\label{eq:ode}
u_j' + \lambda_j (\kappa_1 D_{0+}^\alpha +
\kappa_2 D_{0+}^\beta ) u_j = f_j \quad t\in (0,T),
\end{equation}
subject to the initial value condition $ u_j(0) = \langle u_0,
\phi_j\rangle_{\Omega} $. An elementary computation yields
\begin{equation*}
  u_j + \lambda_j (\kappa_1 I_{0+}^{1-\alpha} +
  \kappa_2 I_{0+}^{1-\beta} ) u_j = I_{0+}f_j+u_j(0)
  \quad \text{ in }  (0,T).
\end{equation*}
Suppose that $ w_j $ satisfies
\[
w_j + \lambda_j (\kappa_1 I_{0+}^{1-\alpha} +
\kappa_2 I_{0+}^{1-\beta} ) w_j = 1  \quad \text{ in } (0,T).
\]
Using the famous Picard iterative process gives
\begin{equation*}
  w_j(t) =  \sum_{r=0}^{\infty}\sum_{p+q=r}\binom{r}{p}
  \frac{
    (-\lambda_j \kappa_1)^p (-\lambda_j \kappa_2)^q\,
    t^{p(1-\alpha)+q(1-\beta)}
  }{
    \Gamma\big(1+p(1-\alpha) + q(1-\beta)\big)
  }.
\end{equation*}
It is easy to verify that
\begin{equation}\label{eq:uj}
  u_j(t) = u_j(0)w_j(t) +\int_{0}^{t}w_j(t-s) f_j(s)\,\mathrm{d}s,
  \quad 0 < t < T,
\end{equation}
which indicates that the singularity part of $ u_j $ belongs to
\[
  S := \text{span}\left\{
    t^{p(1-\alpha) + q(1-\beta)}:\
    p+q > 0,\ p,q \in \mathbb N
  \right\},
\]
provided $ f_j $ is sufficiently regular. Therefore, since
\[
w' \in L_\beta^2(0,T) \quad \text{for all $ w \in S $,}
\]
the assumption $ u' \in L_\beta^2(0,T;H_0^1(\Omega)) $ is
reasonable. We note that the relation \cref{eq:uj} has been applied to ordinary differential equations with multi-term fractional derivatives (\cite{Hadid1996An,Luchko1999AN}).
\section{Proofs}
\label{sec:proof}
\subsection{Auxiliary Results}
Let us first summarize some standard results.
\begin{lem}[\cite{Ciarlet2002,Brenner2008}]
  \label{lem:R_h}
   If $ v \in H_0^1(\Omega) \cap H^{n+1}(\Omega) $, then
  \[
    \nm{(I-R_h)v}_\Omega + h \nm{(I-R_h)v}_{ H_0^1(\Omega) }
    \lesssim h^{n+1} \nm{v}_{ H^{n+1}(\Omega) }.
  \]
  If $ v \in H_0^1(\Omega) $ and $ w \in H^1(\Omega) $, then
  \[
    \dual{(I-R_h)v, w}_\Omega \lesssim h^{\min\{2,n\}}
    \nm{(I-R_h)v}_{H_0^1(\Omega)} \nm{w}_{H_0^1(\Omega)}.
  \]
  If $ v \in H^{m+1}(I_j) $ with $ 1 \leqslant j \leqslant J $, then
  \[
    \nm{(I-Q_\tau)v}_{L^2(I_j)} + \tau_j \nm{(I-Q_\tau)v}_{H^1(I_j)}
    \lesssim \tau_j^{m+1} \nm{v}_{H^{m+1}(I_j)}.
  \]
\end{lem}

\begin{lem}[\cite{Tartar2007}]
  \label{lem:equvi-frac}
  If $ v \in H^\gamma(\mathbb R) $ with $ 0 < \gamma < 1 $, then
  \[
C_\gamma \snm{v}_{H^\gamma(\mathbb R)} \lesssim
    \left(
       \int_\mathbb R \int_\mathbb R
      \frac{\snm{v(s)-v(t)}^2}{\snm{s-t}^{1+2\gamma}}
      \, \mathrm{d}s \, \mathrm{d}t
    \right)^\frac12 \lesssim C_\gamma \snm{v}_{H^\gamma(\mathbb R)}.
  \]
\end{lem}

\begin{lem}[\cite{Samko1993, Diethelm2010}]
  \label{lem:basic-frac}
  The following properties hold:
  \begin{itemize}
    \item If $ 0 < \gamma, \delta < \infty $, then
      \[
        I_{0+}^\gamma I_{0+}^\delta = I_{0+}^{\gamma+\delta}, \quad
        I_{T-}^\gamma I_{T-}^\delta = I_{T-}^{\gamma+\delta}.
      \]
    \item If $ 0 < \gamma< \infty $ and $ u,v \in L^2(0,T) $, then
      \[
        \dual{ I_{0+}^\gamma u, v }_{(0,T)} =
        \dual{ u, I_{T-}^\gamma v }_{(0,T)}.
      \]
  \end{itemize}
\end{lem}

\begin{lem}
  \label{lem:basic}
  If $ v, w \in H^{\gamma/2}(0,T) $ and $ D_{0+}^\gamma v \in L^1(0,T) $ with $
  0 < \gamma < 1 $, then
  \begin{align}
    \dual{D_{0+}^\gamma v, v}_{(0,t)} &
    = \cos(\gamma\pi/2) \snm{v}_{H^{\gamma/2}(0,t)}^2,
    \label{eq:basic-1} \\
    \dual{D_{0+}^\gamma v, w}_{(0,t)} &
    \leqslant  \cos(\gamma\pi/2) \snm{v}_{H^{\gamma/2}(0,t)} \snm{w}_{H^{\gamma/2}(0,t)},
    \label{eq:basic-2}
  \end{align}
  for all $ 0 < t \leqslant T $.
\end{lem}
\noindent The proof of the above lemma is contained in \cite[Lemmas~2.2, 2.4 and
2.9]{Ervin2006}.

The purpose of the rest of this subsection is to prove the following three
lemmas.

\begin{lem}
  \label{lem:vw}
  For $ 0 < t \leqslant T $, it holds that
  \begin{equation}
    \int_0^t \snm{(vw)(s)} \, \mathrm{d}s \lesssim
    \nm{v}_{L_\beta^2(0,t)}
    \snm{w}_{H^{\beta/2}(0,t)}
  \end{equation}
  for all $ v \in L_\beta^2(0,T) $ and $ w \in H^{\beta/2}(0,T) $.
\end{lem}

    \begin{lem}
    \label{lem:jm}
    Let $ 0\leqslant a < b < \infty $ and $ 0 < \gamma < 1 $. If $ v' \in L^2_{\delta}(a,b) $ with $ 0\leqslant \delta < 1 $ and $ v(b) = 0 $, then
    \begin{align}
        \int_a^b v^2(t) (t-a)^{-\gamma} \, \mathrm{d}t
        \leqslant
        \frac{b^{-\delta}}{(1-\delta)(1-\gamma)} (b-a)^{2-\gamma}
        \nm{v'}_{L_\delta^2(a,b)}^2,
        \label{eq:jm-1} \\
        \int_a^b v^2(t) (b-t)^{-\gamma} \, \mathrm{d}t \leqslant
        \frac{b^{-\delta}}{(1-\delta)(1-\gamma)} (b-a)^{2-\gamma}
        \nm{v'}_{L_\delta^2(a,b)}^2,
        \label{eq:jm-2} \\
        \int_a^b \, \mathrm{d}t \int_a^b
        \snm{v(s)-v(t)}^2 \snm{s-t}^{-1-\gamma} \, \mathrm{d}s
        \leqslant \frac{8b^{-\delta}}{1-\delta}(b-a)^{2-\gamma}
        \nm{v'}_{L_\delta^2(a,b)}^2,
        \label{eq:jm-3}
    \end{align}
    where $ B(\cdot,\cdot) $ is the Beta function.
\end{lem}

\begin{lem}
  \label{lem:core}
  For $ 0 < \gamma < 1 $, if $ v \in H^{\gamma/2}(0,T) $ and $ v' \in L^1(0,T)
  $, then
  \begin{equation}\label{eq:core}
    \snm{(I-Q_\tau)v}^2_{H^{\gamma/2}(0,t_j)} \lesssim
    C_\gamma \sum_{i=1}^j \tau_i^{2-\gamma}
    \inf_{0\leqslant\delta<1}\frac{t_i^{-\delta}}{1-\delta}
    \nm{(v-Q_\tau v)'}^2_{L^2_\delta(I_i)}
  \end{equation}
  for each $ 1 \leqslant j \leqslant J $.
\end{lem}

\begin{lem}
  \label{lem:graded}
  If
  \[
    v(t) := t^r, \quad 0\leqslant t \leqslant T,
  \]
  with $ 0 < r \leqslant m +
  1/2 $, then
   \begin{align}
    \label{eq:graded}
      &
      \sum_{i=1}^j \tau_i^{2-\gamma} \inf_{0 \leqslant \delta < 1}
      \frac{t_i^{-\delta}}{1-\delta} \nm{(v-Q_\tau v)'}_{L_\delta^2(I_i)}^2 \nonumber \\
      \lesssim{}& C_{\gamma,\sigma,r}
      \begin{cases}
      T^{(1-\sigma)(2r+1-\gamma)}\tau^{\sigma(2r+1-\gamma)}&
      \text{ if } \sigma < \frac{2m+2-\gamma}{2r+1-\gamma}, \\
      (1+\ln(j) )  T^{2r-1-2m} \tau^{2m+2-\gamma}&
      \text{ if } \sigma = \frac{2m+2-\gamma}{2r+1-\gamma}, \\
      T^{2r-1-2m}\tau^{2m+2-\gamma}
      &
      \text{ if } \sigma > \frac{2m+2-\gamma}{2r+1-\gamma},
      \end{cases}
    \end{align}
  for each $ 1 \leqslant j \leqslant J $.
\end{lem}

\medskip\noindent
{\bf Proof of \cref{lem:vw}.}
  By \cref{lem:equvi-frac}, extending $ w $ to $ \mathbb R \setminus (0,t) $ by
  zero gives
  \[
    \int_0^t \, \mathrm{d}s \int_{-\infty}^0
    w^2(s) (s-\tau)^{-1-\beta} \, \mathrm{d}\tau \lesssim \snm{w}_{H^{\beta/2}(\mathbb R)}^2,
  \]
  which indicates
  \[
    \int_0^t s^{-\beta} w^2(s) \, \mathrm{d}s
    \lesssim \snm{w}_{H^{\beta/2}(\mathbb R)}^2.
  \]
  Therefore, the Cauchy-Schwarz inequality yields
  \begin{align*}
    \int_0^t \snm{(vw)(s)} \, \mathrm{d}s & \leqslant
    \left(
      \int_0^t s^\beta v^2(s) \, \mathrm{d}s
    \right)^\frac12
    \left(
      \int_0^t s^{-\beta} w^2(s) \, \mathrm{d}s
    \right)^\frac12 \\
    & \lesssim  \nm{v}_{L_\beta^2(0,t)}
    \snm{w}_{H^{\beta/2}(\mathbb R)},
  \end{align*}
  and hence the fact $ \snm{w}_{H^{\beta/2}(\mathbb R)} = \snm{w}_{H^{\beta/2}(0,t)} $
  proves the lemma.
\hfill\ensuremath{\blacksquare}

\medskip\noindent
{\bf Proof of \cref{lem:jm}.}
The proof below shall be brief, since the techniques used are standard (see
Minkowski's integral inequality and Hardy's inequality). For $ a < t < b $, a
simple computing gives
\begin{align*}
    {} &
    \snm{v(t)} \leqslant \int_t^b \snm{v'(s)} \, \mathrm{d}s \leqslant
    \left(
    \int_t^b s^{-\delta} \, \mathrm{d}s
    \right)^\frac12
    \left(
    \int_t^b s^\delta \snm{v'(s)}^2 \, \mathrm{d}s
    \right)^\frac12 \\
    \leqslant{} &
    \sqrt{
        \frac{b^{1-\delta} - t^{1-\delta}}{1-\delta}
    }
    \nm{v'}_{L_\delta^2(a,b)} \leqslant
    \sqrt{
        \frac{b^{-\delta}(b-a)}{1-\delta}
    } \nm{v'}_{L_\delta^2(a,b)},
\end{align*}
so that we obtain
\begin{align*}
    {} &
    \int_a^b v^2(t) (t-a)^{-\gamma} \, \mathrm{d}t
    \leqslant   \frac{c(b-a)}{1-\delta}
    \int_a^b  (t-a)^{-\gamma} \, \mathrm{d}t \,
    \nm{v'}_{L_\delta^2(a,b)}^2 \\
    = {} &
    \frac{b^{-\delta}(b- a)^{2-\gamma}}{(1-\delta)(1-\gamma)}
    \nm{v'}_{L_\delta^2(a,b)}^2,
\end{align*}
namely the estimate \cref{eq:jm-1}. Similarly, we have
\begin{align*}
    {} &
    \int_a^b v^2(t) (b-t)^{-\gamma} \, \mathrm{d}t
    \leqslant\frac{b^{-\delta}(b-a)}{1-\delta}
    \int_a^b  (b-t)^{-\gamma}  \, \mathrm{d}t \,
    \nm{v'}_{L_\delta^2(a,b)}^2 \\
    = {} &
    \frac{b^{-\delta}(b- a)^{2-\gamma}}{(1-\delta)(1-\gamma)}
    \nm{v'}_{L_\delta^2(a,b)}^2,
\end{align*}
namely the estimate \cref{eq:jm-2}.
Finally, let us prove \cref{eq:jm-3}. Since
\begin{align*}
    {} &
    \int_a^b \, \mathrm{d}t
    \int_a^b \snm{v(s)-v(t)}^2 \snm{s-t}^{-1-\gamma} \, \mathrm{d}s \\
    ={} &
    2\int_a^b \, \mathrm{d}t
    \int_t^b \snm{\int_t^s v'(\tau) \, \mathrm{d}\tau}^2 (s-t)^{-1-\gamma}
    \, \mathrm{d}s \\
    ={} &
    2\int_a^b \, \mathrm{d}t
    \int_t^b \snm{\int_0^1 v'(t+\theta(s-t)) \, \mathrm{d}\theta}^2
    (s-t)^{1-\gamma} \, \mathrm{d}s \\
    \leqslant{} &
    2(b-a)^{1-\gamma} \int_{a}^{b}
    \left(
    \int_0^1 \sqrt{
        \int_t^b \snm{v'(t+\theta(s-t))}^2 \, \mathrm{d}s
    } \, \mathrm{d}\theta
    \right)^2 \, \mathrm{d} t \\
    ={} &
    2(b-a)^{1-\gamma} \int_a^b
    \left(
    \int_0^1 \sqrt{
        \int_t^{t+\theta(b-t)}
        \snm{v'(\eta)}^2 \theta^{-1} \,\mathrm{d}\eta
    } \,\mathrm{d}\theta
    \right)^2 \, \mathrm{d}t,
\end{align*}
the inequality \cref{eq:jm-3} is a direct consequence of
\begin{align*}
    {} &
    \int_a^b \left(
    \int_0^1 \sqrt{
        \int_t^{t+\theta(b-t)}
        \snm{v'(\eta)}^2 \theta^{-1} \, \mathrm{d}\eta
    } \, \mathrm{d}\theta
    \right)^2 \, \mathrm{d}t \\
    \leqslant{} &
    \int_a^b \left(
    \int_0^1 \sqrt{
        \int_t^{t+\theta(b-t)}
        (\eta/t)^\delta \snm{v'(\eta)}^2 \theta^{-1} \, \mathrm{d}\eta
    } \,\mathrm{d}\theta
    \right)^2 \, \mathrm{d}t \\
    ={} &
    \int_a^b t^{-\delta} \left(
    \int_0^1 \theta^{-1/2} \sqrt{
        \int_t^{t+\theta(b-t)}
        \eta^\delta \snm{v'(\eta)}^2 \, \mathrm{d}\eta
    } \,\mathrm{d}\theta
    \right)^2 \, \mathrm{d}t \\
    \leqslant{} &
    \int_a^b t^{-\delta} \left(
    \int_0^1 \theta^{-1/2}
    \,\mathrm{d}\theta
    \right)^2 \, \mathrm{d}t \, \nm{v'}_{L_\delta^2(a,b)}^2 \\
    \leqslant{} &
    \frac{4b^{-\delta}(b-a)}{1-\delta}
    \nm{v'}_{L_\delta^2(a,b)}^2.
\end{align*}
This lemma is thus proved.
\hfill\ensuremath{\blacksquare}

\medskip\noindent
{\bf Proof of \cref{lem:core}.}
By \cref{lem:equvi-frac} we only need to prove
  \begin{equation}
    \label{eq:I123}
    \begin{aligned}
       \mathbb I_1 + \mathbb I_2 + \mathbb I_3
    \lesssim C_\gamma \sum_{i=1}^j \tau_i^{2-\gamma}\inf_{0\leqslant\delta<1}\frac{t_i^{-\delta}}{1-\delta}\nm{(v-Q_\tau v)'}^2_{L^2_\delta(I_i)}.
    \end{aligned}
  \end{equation}
where
  \begin{align*}
    \mathbb I_1 &= \sum_{i=1}^j
    \int_{t_{i-1}}^{t_i} \, \mathrm{d}t
    \int_{t_{i-1}}^{t_i} \snm{g(t)-g(s)}^2\snm{t-s}^{-1-\gamma} \, \mathrm{d}s, \\
    \mathbb I_2 &= \sum_{i=1}^j \sum_{l=i+1}^j
    \int_{t_{i-1}}^{t_i} \, \mathrm{d}t
    \int_{t_{l-1}}^{t_l} \snm{g(t)-g(s)}^2\snm{t-s}^{-1-\gamma} \mathrm{d}s, \\
    \mathbb I_3 &= \int_0^{t_j} \snm{g(t)}^2
    \left(
      \int_{t_j}^\infty (s-t)^{-1-\gamma} \,\mathrm{d}s +
      \int_{-\infty}^0 (t-s)^{-1-\gamma} \,\mathrm{d}s
    \right)\, \mathrm{d}t.
  \end{align*}

  A straightforward calculation gives
  \begin{align*}
    {} &
    \sum_{i=1}^j \sum_{l=i+1}^j
    \int_{t_{i-1}}^{t_i} \, \mathrm{d}t \int_{t_{l-1}}^{t_l}
    g^2(t)\snm{t-s}^{-1-\gamma} \, \mathrm{d}s \\
    ={} &
    \frac1{\gamma} \sum_{i=1}^j \sum_{l=i+1}^j \int_{t_{i-1}}^{t_i}
    g^2(t) \left(
      (t_{l-1}-t)^{-\gamma} - (t_l-t)^{-\gamma}
    \right)
    \, \mathrm{d}t  \\
    \leqslant{} &
    \frac1{\gamma} \sum_{i=1}^{j-1} \int_{t_{i-1}}^{t_i}\
    g^2(t) (t_i-t)^{-\gamma}
    \, \mathrm{d}t
  \end{align*}
  and
  \begin{align*}
    {} &
    \sum_{i=1}^j \sum_{l=i+1}^j
    \int_{t_{i-1}}^{t_i} \, \mathrm{d}t \int_{t_{l-1}}^{t_l}
    g^2(s)\snm{t-s}^{-1-\gamma} \, \mathrm{d}s \\
    ={} &
    \frac1\gamma \sum_{i=1}^j \sum_{l=i+1}^j \int_{t_{l-1}}^{t_l} g^2(s)
    \left(
      (s-t_i)^{-\gamma} - (s-t_{i-1})^{-\gamma}
    \right) \, \mathrm{d}s \\
    \leqslant{} &
    \frac1{\gamma} \sum_{l=2}^j \int_{t_{l-1}}^{t_l}
    g^2(s) (s-t_{l-1})^{-\gamma}
    \, \mathrm{d}s.
  \end{align*}
  It follows
  \[
    \mathbb I_2 \leqslant
    \frac2\gamma \sum_{i=1}^j
    \int_{t_{i-1}}^{t_i} g^2(t) \big(
      (t_i-t)^{-\gamma} + (t-t_{i-1})^{-\gamma}
    \big) \, \mathrm{d}t.
  \]
  Therefore, since it is evident that
  \begin{align*}
    \mathbb I_3 \leqslant
    \frac1\gamma \sum_{i=1}^j \int_{t_{i-1}}^{t_i}
    g^2(t) \big((t_i-t)^{-\gamma} + (t-t_{i-1})^{-\gamma}\big) \, \mathrm{d}t,
  \end{align*}
  using \cref{lem:jm} yields
  \begin{align*}
    \mathbb I_2 + \mathbb I_3 &
    \leqslant \frac3\gamma \sum_{i=1}^j
    \int_{t_{i-1}}^{t_i} g^2(t)
    \big( (t_i-t)^{-\gamma} + (t-t_{i-1})^{-\gamma} \big)
    \, \mathrm{d}t \\
    & \lesssim
 C_\gamma\sum_{i=1}^j \tau_i^{2-\gamma}\inf_{0\leqslant\delta<1}\frac{t_i^{-\delta}}{1-\delta}\nm{(v-Q_\tau v)'}^2_{L^2_\delta(I_i)}.
  \end{align*}
  As using \cref{lem:jm} also yields
  \[\begin{aligned}
    \mathbb I_1 & \lesssim C_\gamma\sum_{i=1}^j \tau_i^{2-\gamma}\inf_{0\leqslant\delta<1}\frac{t_i^{-\delta}}{1-\delta}\nm{(v-Q_\tau v)'}^2_{L^2_\delta(I_i)},
    \end{aligned}
  \]
  we readily obtain \cref{eq:I123} and thus complete the proof of
  \cref{lem:core}.
\hfill\ensuremath{\blacksquare}

\medskip\noindent
{\bf Proof of \cref{lem:graded}.}
  Setting
  \[
    \delta_0 :=
    \begin{cases}
      1-r & \text{ if } 0 < r < 1/2, \\
      1/2 & \text{ if } r \geqslant 1/2,
    \end{cases}
  \]
  by a standard scaling argument we obtain
  \[
    \nm{(Q_\tau v)'}_{L_{\delta_0}^2(I_1)} \lesssim
    C_r \nm{v'}_{L_{\delta_0}^2(I_1)} \lesssim
    C_r \tau_1^{(2r+\delta_0-1)/2}
  \]
  and hence
  \[
    \inf_{0 \leqslant \delta < 1} \frac{t_1^{-\delta}}{1-\delta}
    \nm{(v-Q_\tau v)'}_{L_\delta^2(I_1)}^2 \leqslant
    \frac{\tau_1^{-\delta_0}}{1-\delta_0} \nm{(v-Q_\tau v)'}_{L_{\delta_0}^2(I_1)}^2
    \lesssim C_r \tau_1^{2r-1}.
  \]
  Therefore, \cref{lem:R_h} implies
  \begin{align}
    &
    \sum_{i=1}^j \tau_i^{2-\gamma} \inf_{0 \leqslant \delta < 1}
    \frac{t_i^{-\delta}}{1-\delta} \nm{(v-Q_\tau v)'}_{L_\delta^2(I_i)}^2 \nonumber \\
    \lesssim{} &
    C_r \tau_1^{2r+1-\gamma} + \sum_{i=2}^j \tau_i^{2m+2-\gamma}
    \nm{v^{(m+1)}}_{L^2(I_i)}^2\label{eq:est_inf}.
  \end{align}
  Since a simple computing yields
  \[
    \tau_i < \frac{2^{\sigma-1}}\sigma J^{-1} T^{1/\sigma}
    t_{i-1}^{1-1/\sigma}, \quad 2 \leqslant i \leqslant J,
  \]
  it follows
  \begin{align}
    {} &
    \sum_{i=2}^j \tau_i^{2m+2-\gamma} \nm{v^{(m+1)}}_{L^2(I_i)}^2 \nonumber \\
    \leqslant{} &C_{\sigma,r}
    \left(
      J^{-1} T^{1/\sigma}
    \right)^{2m+2-\gamma}
    \int_{t_1}^{t_j} t^{2r-\gamma-(2m+2-\gamma)/\sigma}
    \, \mathrm{d}t \nonumber \\
    \leqslant{} & C_{\gamma,\sigma,r}
        \begin{cases}
          T^{2r+1-\gamma}J^{-\sigma(2r+1-\gamma)}&
          \text{ if } \sigma < \frac{2m+2-\gamma}{2r+1-\gamma}, \\
          \ln\left(t_j/t_1 \right) T^{2r+1-\gamma}J^{-2m-2+\gamma}&
          \text{ if } \sigma = \frac{2m+2-\gamma}{2r+1-\gamma}, \\
          t_j^{2r+1-\gamma-(2m+2-\gamma)/\sigma}T^{(2m+2-\gamma)/\sigma}J^{-2m-2+\gamma}
          &
          \text{ if } \sigma > \frac{2-2m-\gamma}{1+2r-\gamma}.
        \end{cases}\label{eq:est}
  \end{align}
  Therefore, by \cref{eq:est_inf,eq:est} and the fact $ T/J < \tau $, a direct
  computation yields \cref{eq:graded} and thus concludes the proof of
  \cref{lem:graded}.
\hfill\ensuremath{\blacksquare}

\subsection{Proofs of \cref{thm:stability,thm:conv,coro:conv,coro:graded}}
Since the proofs of \cref{coro:conv,coro:graded} are
straightforward by \cref{thm:conv,lem:R_h,lem:graded}, this subsection only
proves \cref{thm:stability,thm:conv}. For $ 1 \leqslant j \leqslant J $, set
$ \Omega_{t_j} := \Omega \times (0,t_j) $ and define
\[
  s_j(V,W) := \sum_{i=0}^{j-1} \dual{\jmp{V_j}, W_j^{+}}_{\Omega},
  \quad \forall\, V, W \in W_{h,\tau}.
\]
\begin{lem}
  \label{lem:sj}
  If $ V \in W_{h,\tau} $ and $ v \in L^2(\Omega) $, then
  \[
    \frac12 \left(
      \nm{V_j}_\Omega^2 + \nm{V_0^{+}}_\Omega^2
    \right) \leqslant
    \dual{V',V}_{\Omega_{t_j}} + s_j(V,V),
  \]
  \[
    \frac12 \left(
      \nm{V_j}_\Omega^2 - \nm{v}_\Omega^2
    \right) \leqslant
    \dual{V',V}_{\Omega_{t_j}} + s_j(V,V) - \dual{v,V_0^{+}}_{\Omega},
  \]
  for all $ 1 \leqslant j \leqslant J $.
\end{lem}
\noindent The above lemma is contained in the proof of
\cite[Theorem~12.1]{Thomee2006}.

\medskip\noindent
{\bf Proof of \cref{thm:stability}.} Since the stability result  \cref{eq:stab} indicates the unique existence of $U$, it suffices to prove the former.
  By \cref{lem:basic,lem:sj}, inserting $ V = U \chi_{(0,t_j)} $ into
  \cref{eq:U} yields
  \begin{align*}
    & \frac12 \nm{U_j}_{L^2(\Omega)}^2 +
    \kappa_1 \snm{U}_{H^{\alpha/2}(0,t_j;H_0^1(\Omega))}^2 +
    \kappa_2 \snm{U}_{H^{\beta/2}(0,t_j;H_0^1(\Omega))}^2 \\
    \leqslant{} &
    \frac12 \nm{R_hu_0}_{L^2(\Omega)}^2 +
    \dual{f,U}_{L^\infty(0,t_j;H_0^1(\Omega))},
  \end{align*}
  so that \cref{lem:vw} implies
  \begin{align*}
    & \nm{U_j}_{L^2(\Omega)}^2 +
    \kappa_1 \snm{U}_{H^{\alpha/2}(0,t_j;H_0^1(\Omega))}^2 +
    \kappa_2 \snm{U}_{H^{\beta/2}(0,t_j;H_0^1(\Omega))}^2 \\
    \lesssim{} &
    \nm{R_hu_0}_{L^2(\Omega)}^2 +
    \nm{f}_{L_\beta^2(0,t_j;H^{-1}(\Omega))}
    \snm{U}_{H^{\beta/2}(0,t_j;H_0^1(\Omega))}.
  \end{align*}
  Therefore, using the Young's inequality with $ \epsilon $ gives
  \begin{align*}
    &\nm{U_j}_{L^2(\Omega)} +
    \sqrt{\kappa_1} \snm{U}_{H^{\alpha/2}(0,t_j;H_0^1(\Omega))} +
    \sqrt{\kappa_2} \snm{U}_{H^{\beta/2}(0,t_j;H_0^1(\Omega))} \\
    \lesssim{} &
    \nm{R_hu_0}_{L^2(\Omega)} +
    1/\sqrt{\kappa_2} \nm{f}_{L_\beta^2(0,t_j;H^{-1}(\Omega))},
  \end{align*}
  which, together with the estimate
  \[
    \nm{R_hu_0}_{L^2(\Omega)} \lesssim \nm{u_0}_{H_0^1(\Omega)},
  \]
  proves \cref{eq:stab} and thus concludes the proof of \cref{thm:stability}.
\hfill\ensuremath{\blacksquare}

\medskip\noindent
{\bf Proof of \cref{thm:conv}.}
  By integration by parts, using \cref{eq:model} yields
  \[
    \dual{f,\theta}_{L^\infty(0,t_j;H_0^1(\Omega))} =
    \dual{u',\theta}_{\Omega_{t_j}} +
    \dual{
      \left(
        \kappa_1 D_{0+}^\alpha + \kappa_2 D_{0+}^\beta
      \right) \partial_xu,\partial_x\theta
    }_{\Omega_{t_j}}.
  \]
  Moreover, substituting $ V = \theta \chi_{(0,t_j)} $ into \cref{eq:U} yields
  \begin{align*}
    \dual{f,\theta}_{L^\infty(0,t_j;H_0^1(\Omega))} = &
    \dual{U',\theta}_{\Omega_{t_j}} +
    \dual{
      \left(
        \kappa_1 D_{0+}^\alpha + \kappa_2 D_{0+}^\beta
      \right) \partial_xU,\partial_x\theta
    }_{\Omega_{t_j}} \\
    & {} + s_j(U,\theta) - \dual{R_hu_0,\theta_0^{+}}_\Omega.
  \end{align*}
  Consequently, it follows
  \begin{align*}
    0 = & \dual{(u-U)',\theta}_{\Omega_{t_j}} +
    \dual{
      \left(
        \kappa_1 D_{0+}^\alpha + \kappa_2 D_{0+}^\beta
      \right)
      \partial_x(u-U),\partial_x\theta
    }_{\Omega_{t_j}} \\
    & s_j(u-U,\theta) - \dual{(I-R_h)u_0,\theta_0^{+}}_\Omega,
  \end{align*}
  and then a simple calculation gives
  \[
    \dual{\theta',\theta}_{\Omega_{t_j}} +
    \dual{
      \left(
        \kappa_1 D_{0+}^\alpha + \kappa_2 D_{0+}^\beta
      \right) \partial_x\theta,\partial_x\theta
    }_{\Omega_{t_j}} +
    s_j(\theta,\theta)  = \mathbb I_1 + \mathbb I_2 + \mathbb I_3,
  \]
  where $ \rho := (I-Q_\tau R_h) u $ and
  \begin{align*}
    \mathbb I_1 &= \dual{\rho', \theta}_{\Omega_{t_j}} + s_j(\rho,\theta) -
    \dual{(I-R_h)u_0,\theta_0^{+}}_\Omega, \\
    \mathbb I_2 &= \kappa_1 \dual{D_{0+}^\alpha\partial_x\rho,
    \partial_x\theta}_{\Omega_{t_j}}, \\
    \mathbb I_3 &= \kappa_2 \dual{D_{0+}^\beta\partial_x\rho,
    \partial_x\theta}_{\Omega_{t_j}}.
  \end{align*}
  Therefore, \cref{lem:sj} implies
  \[
    \nm{\theta_j}_{L^2(\Omega)}^2 +
    \kappa_1 \snm{\theta}_{H^{\alpha/2}(0,t_j;H_0^1(\Omega))}^2 +
    \kappa_2 \snm{\theta}_{H^{\beta/2}(0,t_j;H_0^1(\Omega))}^2
    \lesssim \mathbb I_1 + \mathbb I_2 + \mathbb I_3.
  \]

  Let us first estimate $ \mathbb I_1 $. By the definition of $ Q_\tau $, using
  integration by parts gives
  \begin{align*}
    &
    \dual{((I-Q_\tau)R_hu)', \theta}_{\Omega_{t_j}} +
    s_j((I-Q_\tau)R_hu, \theta) \\
    ={} &
    \sum_{i=1}^j -\dual{((I-Q_\tau)R_hu)_{i-1}^{+},\theta_{i-1}^{+}}_\Omega +
    \sum_{i=1}^{j-1} \dual{(I-Q_\tau)R_hu)_i^{+}, \theta_i^{+}}_\Omega \\
    {} & {} +
    \dual{((I-Q_\tau)R_hu)_0^{+}n, \theta_0^{+}}_\Omega = 0,
  \end{align*}
  which implies
  \begin{align*}
    \mathbb I_1 &=
    \dual{(u-R_hu)', \theta}_{\Omega_{t_j}} + s_j(u-R_hu,\theta) -
    \dual{(I-R_h)u_0, \theta_0^{+}}_\Omega \\
    &=
    \dual{(I-R_h)u', \theta}_{\Omega_{t_j}}.
  \end{align*}
  Therefore, using \cref{lem:R_h,lem:vw} yields
  \[
    \mathbb I_1 \lesssim \eta_{j,1} \sqrt{\kappa_2}
    \snm{\theta}_{H^{\beta/2}(0,t_j;H_0^1(\Omega))}.
  \]

  Then let us estimate $ \mathbb I_2 $ and $ \mathbb I_3 $. A straightforward
  calculation gives
  \begin{align*}
    \mathbb I_2 &= \kappa_1 \dual{
      D_{0+}^\alpha \partial_x (u-Q_\tau R_hu),
      \partial_x\theta
    }_{\Omega_{t_j}} \\
    &= \kappa_1 \dual{
      D_{0+}^\alpha \partial_x(u-Q_\tau u),
      \partial_x\theta
    }_{\Omega_{t_j}} \\
    &= \kappa_1 \dual{
      D_{0+}^\alpha (I-Q_\tau)\partial_xu,
      \partial_x\theta
    }_{\Omega_{t_j}},
  \end{align*}
  so that \cref{lem:basic,lem:core} imply
  \[
    \mathbb I_2 \lesssim \eta_{j,2}
    \sqrt{\kappa_1} \snm{\theta}_{H^{\alpha/2}(0,t_j;H_0^1(\Omega))}.
  \]
  Analogously, we obtain
  \[
    \mathbb I_3 \lesssim \eta_{j,3}
    \sqrt{\kappa_2} \snm{\theta}_{H^{\beta/2}(0,t_j;H_0^1(\Omega))}.
  \]

  Finally, by the Young's inequality with $ \epsilon $, combining the above
  estimates for $ \mathbb I_1 $, $ \mathbb I_2 $ and $ \mathbb I_3 $ yields
  \cref{eq:conv}. This concludes the proof of \cref{thm:conv}.
\hfill\ensuremath{\blacksquare}

\section{Numerical Experiments}
\label{sec:numer}
This section investigates numerically the temporal accuracy of $ U $. We set $
\alpha = 0.2 $, $ \beta = 0.8 $, $ \kappa_1 = \kappa_2 = 1 $, $ \Omega = (0,1) $
and $ T = 1 $, and let
\[
  u(x,t) := t^r \sin(\pi x), \quad  (x,t) \in \Omega_T
\]
be the solution to problem \cref{eq:model}, where $ r > 0 $ is a constant. To
ensure that the spatial discretization error is negligible compared with the
temporal discretization error, we set $ n = 3 $ and $ h = 1/32 $. Additionally,
define
\begin{align*}
  E_1(U) &:= \max_{1 \leqslant j \leqslant J} \nm{(u-U)(t_j)}_{L^2(\Omega)} ,\\
  E_2(U) &:= \nm{(u-U)(T)}_{L^2(\Omega)}.
\end{align*}

\noindent {\bf Experiment 1.} This experiment investigates the temporal accuracy
of $ U $ under the condition that $ u $ is sufficiently regular and the temporal
grid is equidistant ($ \sigma = 1 $). We set $ r = 4 $ and present the
corresponding numerical results in \cref{tab:ex1}. These numerical results show
that $ E_1(U) = O(\tau^{m+1}) $, which exceeds the theoretical temporal accuracy
$ O(\tau^{m+0.6}) $ indicated by \cref{coro:conv}.

\begin{table}[H]
    \caption{$ r=4 ,~\sigma = 1 $.}
    \label{tab:ex1}
    \begin{tabular}{cccccc}
        \toprule
        \multirow{2}{*}{$J$} &
        \multicolumn{2}{c}{$m=0$} &&
        \multicolumn{2}{c}{$m=1$} \\
        \cmidrule{2-3} \cmidrule{5-6}
                 & $E_1(U)$ & Order && $E_1(U)$ & Order \\
        \midrule
        $ 64  $  &  1.43e-2 & --    &&  3.02e-5 & --    \\
        $ 128 $  &  7.56e-3 & 0.92  &&  7.20e-6 & 2.07  \\
        $ 256 $  &  3.95e-3 & 0.93  &&  1.71e-6 & 2.07  \\
        $ 512 $  &  2.05e-3 & 0.95  &&  4.08e-7 & 2.07  \\
        $ 1024 $ &  1.06e-3 & 0.96  &&  9.69e-8 & 2.07  \\
        \bottomrule
    \end{tabular}
\end{table}

\noindent{\bf Experiment 2.} This experiment investigates the temporal accuracy
of $ U $ under the condition that $ u $ has singularity near $ t={0+} $ and the
temporal grid is also equidistant. The corresponding numerical results are
displayed in \cref{tab:singu-m0,tab:singu-m1}, and they illustrate that $ E_1(U)
= O(\tau^{r+0.1}) $ which agrees with \cref{coro:graded}. The numerical results
also show that the theoretical accuracy
\[
  E_1(U)= O(\tau^{m+1-\beta/2})
\]
indicated by \cref{coro:conv} is optimal with respect to the regularity of $ u
$. Furthermore, \cref{tab:singu-m0,tab:singu-m1} illustrate the following
interesting result:
\[
  E_2(U) = O(\tau^{m+1}).
\]
Therefore, if only $ u(T) $ is concerned, then equidistant temporal grids are
sufficient.

\begin{table}[H]
    \caption{$ m=0 ,~\sigma = 1 $.}
    \label{tab:singu-m0}
    \begin{tabular}{ccccccc}
        \toprule
        $r$ & $J$  & $E_1(U)$ & Order && $E_2(U)$ & Order \\
        \midrule
        \multirow{5}{*}{0.2}
        &$ 16  $   &  2.32e-2 & --    && 5.29e-3  & --    \\
        &$ 32  $   &  1.83e-2 & 0.34  && 2.64e-3  & 1.00  \\
        &$ 64  $   &  1.46e-2 & 0.32  && 1.32e-3  & 1.00  \\
        &$ 128 $   &  1.19e-2 & 0.30  && 6.55e-4  & 1.01  \\
        &$ 256 $   &  9.73e-3 & 0.29  && 3.26e-4  & 1.01  \\
        \midrule
        \multirow{5}{*}{0.5}
        &$ 16  $   &  2.09e-2 & --    && 1.09e-2  & --    \\
        &$ 32  $   &  1.34e-2 & 0.64  && 5.51e-3  & 0.99  \\
        &$ 64  $   &  8.75e-3 & 0.62  && 2.76e-3  & 0.99  \\
        &$ 128 $   &  5.78e-3 & 0.60  && 1.38e-3  & 1.00  \\
        &$ 256 $   &  3.85e-3 & 0.59  && 6.93e-4  & 1.00  \\
        \midrule
        \multirow{5}{*}{0.8}
        &$ 16 $    &  1.59e-2 & --    && 1.55e-2  & --    \\
        &$ 32 $    &  8.41e-3 & 0.92  && 7.89e-3  & 0.97  \\
        &$ 64 $    &  4.45e-3 & 0.92  && 4.01e-3  & 0.98  \\
        &$ 128 $   &  2.37e-3 & 0.91  && 2.03e-3  & 0.98  \\
        &$ 256 $   &  1.27e-3 & 0.90  && 1.02e-3  & 0.99  \\
        \bottomrule
    \end{tabular}
\end{table}

\begin{table}[H]
    \caption{$ m=1,~\sigma = 1 $.}
    \label{tab:singu-m1}
    \begin{tabular}{ccccccc}
        \toprule
        $r$ & $J$  & $E_1(U)$ & Order && $E_2(U)$ & Order \\
        \midrule
        \multirow{5}{*}{0.5}
        &$ 16  $   &  1.89e-3 & --    && 1.42e-5  & --    \\
        &$ 32  $   &  1.23e-3 & 0.62  && 3.10e-6  & 2.20  \\
        &$ 64  $   &  8.10e-4 & 0.60  && 6.97e-7  & 2.15  \\
        &$ 128 $   &  5.38e-4 & 0.59  && 1.60e-7  & 2.12  \\
        &$ 256 $   &  3.59e-4 & 0.58  && 3.70e-8  & 2.11  \\
        \midrule
        \multirow{5}{*}{0.8}
        &$ 16 $    &  4.08e-4 & --    && 8.36e-6  & --    \\
        &$ 32 $    &  2.16e-4 & 0.92  && 1.86e-6  & 2.17  \\
        &$ 64 $    &  1.15e-4 & 0.90  && 4.26e-7  & 2.13  \\
        &$ 128 $   & 6.23e-5  & 0.89  && 9.90e-8  & 2.11  \\
        &$ 256 $   &  3.38e-5 & 0.88  && 2.30e-8  & 2.10  \\
        \midrule
        \multirow{5}{*}{1.5}
        &$ 16 $    &  1.71e-4 & --    && 3.65e-5  & --    \\
        &$ 32 $    &  5.57e-5 & 1.62  && 8.36e-6  & 2.12  \\
        &$ 64 $    &  1.84e-5 & 1.60  && 1.95e-6  & 2.10  \\
        &$ 128 $   &  6.13e-6 & 1.59  && 4.57e-7  & 2.09  \\
        &$ 256 $   &  2.05e-6 & 1.58  && 1.08e-7  & 2.08  \\
        \bottomrule
    \end{tabular}
\end{table}

\noindent{\bf Experiment 3.} This experiment investigates the temporal accuracy
of $ U $ under the condition that $ u $ has singularity near $ t={0+} $ and the
temporal grid is graded with different parameter $ \sigma > 1 $. We consider $ r
= 0.2 $ and $ r = 0.4 ,$ and list the corresponding numerical results in
\cref{tab:ex3r1m0,tab:ex3r1m1,tab:ex3r2m0,tab:ex3r2m1}. For $ 1 < \sigma
\leqslant \sigma^* $, the numerical results show that $ E_1(U) =
O(\tau^{\sigma(r+0.1)}) $, which agrees with \cref{coro:graded}. Moreover, in
the case of $ \sigma = \sigma^{**} $, the temporal accuracy $ E_1(U) =
O(\tau^{m+1}) $ is observed. Here, we recall that $ \sigma^* $ and $ \sigma^{**}
$ are defined by \cref{eq:s1,eq:s2}, respectively.

\begin{table}[H]
    \caption{$ r=0.2 $, $ m=0 $.}
    \label{tab:ex3r1m0}
    \begin{tabular}{ccccccc}
        \toprule
        $\sigma\backslash J$ &&16&32& 64 & 128 & 256 \\

        \midrule
        \multirow{2}{*}{1.5}
        & $E_1(U)$ & 1.46e-02 & 1.07e-02 & 8.02e-03 & 6.04e-03 & 4.55e-03 \\
        & Order    & --       & 0.45     & 0.42     & 0.41     & 0.41     \\

        \midrule
        \multirow{2}{*}{2($ \sigma^* $)}
        & $E_1(U)$ & 1.05e-02 & 7.04e-03 & 4.82e-03 & 3.31e-03 & 2.26e-03 \\
        & Order    & --       & 0.58     & 0.55     & 0.54     & 0.55     \\

        \midrule
        \multirow{2}{*}{$\dfrac{10}{3}$($\sigma^{**}$)}
        & $E_1(U)$ & 1.11e-02 & 5.87e-03 & 3.05e-03 & 1.57e-03 & 8.05e-04 \\
        & Order    & --       & 0.92     & 0.95     & 0.96     & 0.96     \\

        \bottomrule
    \end{tabular}
\end{table}

\begin{table}[H]
    \caption{$ r=0.2 $, $ m=1 $.}
    \label{tab:ex3r1m1}
    \begin{tabular}{ccccccc}
        \toprule
        $\sigma\backslash J$ &&4&8&16&32 & 64  \\

        \midrule
        \multirow{2}{*}{3}
        & $E_1(U)$ & 2.44e-03 & 1.35e-03 & 7.40e-04 & 3.94e-04 & 2.01e-04 \\
        & Order    & --       & 0.86     & 0.86     & 0.91     & 0.97     \\

        \midrule
        \multirow{2}{*}{$\dfrac{16}{3}$($\sigma^*$)}
        & $E_1(U)$ & 3.73e-03 & 1.01e-03 & 3.35e-04 & 1.10e-04 & 3.30e-05 \\
        & Order    & --       & 1.89     & 1.59     & 1.61     & 1.73     \\

        \midrule
        \multirow{2}{*}{$\dfrac{20}{3}$($\sigma^{**}$)}
        & $E_1(U)$ & 5.13e-03 & 1.51e-03 & 3.81e-04 & 9.11e-05 & 2.33e-05 \\
        & Order    & --       & 1.76     & 1.99     & 2.06     & 1.97     \\

        \bottomrule
    \end{tabular}
\end{table}

\begin{table}[H]
    \caption{$ r=0.4 $, $ m=0 $.}
    \label{tab:ex3r2m0}
    \begin{tabular}{ccccccc}
        \toprule
        $\sigma\backslash J$ &&16&32& 64 & 128 & 256 \\

        \midrule
        \multirow{2}{*}{1.1}

        & $E_1(U)$ & 2.01e-02 & 1.34e-02 & 9.11e-03 & 6.27e-03 & 4.35e-03 \\
        & Order    & --       & 0.58     & 0.56     & 0.54     & 0.53  \\

        \midrule
        \multirow{2}{*}{1.2($ \sigma^* $)}
        & $E_1(U)$ & 1.73e-02 & 1.12e-02 & 7.42e-03 & 4.97e-03 & 3.34e-03 \\
        & Order    & --       & 0.62     & 0.60     & 0.58     & 0.57  \\

        \midrule
        \multirow{2}{*}{2($\sigma^{**}$)}

        & $E_1(U)$ & 1.41e-02 & 7.31e-03 & 3.76e-03 & 1.92e-03 & 9.85e-04 \\
        & Order    & --       & 0.94     & 0.96     & 0.97     & 0.97 \\

        \bottomrule
    \end{tabular}
\end{table}

\begin{table}[H]
    \caption{$ r=0.4 $, $ m=1 $.}
    \label{tab:ex3r2m1}
    \begin{tabular}{ccccccc}
        \toprule
        $\sigma\backslash J$ &&4&8&16&32 & 64  \\

        \midrule
        \multirow{2}{*}{2}
        & $E_1(U)$ & 2.64e-03 & 1.29e-03 & 6.59e-04 & 3.37e-04 & 1.71e-04 \\
        & Order    & --       & 1.03     & 0.98     & 0.97     & 0.98  \\

        \midrule
        \multirow{2}{*}{3.2($\sigma^*$)}

        & $E_1(U)$ & 2.27e-03 & 6.69e-04 & 2.29e-04 & 7.73e-05 & 2.53e-05\\
        & Order    & --       & 1.76     & 1.55     & 1.57     & 1.61  \\

        \midrule
        \multirow{2}{*}{4($\sigma^{**}$)}
        & $E_1(U)$ & 3.28e-03 & 8.43e-04 & 2.02e-04 & 4.81e-05 & 1.26e-05\\
        & Order    & --       & 1.96     & 2.06     & 2.07     & 1.93  \\

        \bottomrule
    \end{tabular}
\end{table}

\section{Conclusions}
\label{sec:concluding}
This paper analyzes a time-stepping discontinuous Galerkin method for the modified
anomalous subdiffusion problem. We establish the stability of this method and
prove that the temporal accuracy is $ O(\tau^{m+1-\beta/2}) $, and the numerical
results confirm that this accuracy is optimal with respect to the regularity of
$ u $. Furthermore, if $ u $ has singularity near $ t = {0+} $, we prove that
employing graded grids in the temporal discretization can improve the temporal
accuracy to $ O(\tau^{m+1-\beta/2}) $, which is also verified by the numerical
results.

However, further investigations are still needed.
\begin{itemize}
    \item The numerical results illustrate that if $ u $ is sufficiently regular,
    then 
    $$\max_{1 \leqslant j \leqslant J} \nm{(u-U)(t_j)}_{L^2(\Omega)}= O(\tau^{m+1}) .$$
    \item Although $ u $ has singularity near $ t = {0+} $, the numerical results
    show that
    \[
    \nm{(u-U)(T)}_{L^2(\Omega)} = O(\tau^{m+1}).
    \]
    \item The numerical results also illustrate that if $ u $ has singularity near
    $ t = {0+} $, then adopting graded grids in the temporal discretization can
    improve the temporal accuracy to $ O(\tau^{m+1}) $.
\end{itemize}


\end{document}